\documentclass[ctagsplt,12pt]{amsart}
\RequirePackage{newlfont}
\usepackage{amscd}
\usepackage{amsmath}
\usepackage{amsfonts}
\usepackage{amssymb}
\usepackage{enumerate}
\usepackage{graphicx}
\usepackage{latexsym}
\usepackage{color}
\usepackage{bbm}

\usepackage{relsize} % To use \mathlarger for larger sizes

\numberwithin{equation}{section}

\begin{document}

	\title[Limits]{$p(x)$-Stability of the Dirichlet problem for Poisson's equation with variable exponents}
	
	\author[ B. Djafari Rouhani, O. M\'{e}ndez ]{Behzad Djafari Rouhani and Osvaldo M\'{e}ndez }

	\address{Osvaldo M\'{e}ndez\\Department of Mathematical Sciences, The University of Texas at El Paso, El Paso, TX 79968, USA}
	\email{osmendez@utep.edu}
	\address{Behzad Djafari Rouhani\\ Department of Mathematical Sciences, University of Texas at El Paso}
	\email{behzad@utep.edu}
	
	
	\begin{abstract}
		It is shown that if the sequence $(p_j(x))$ increases uniformly to $p(x)$ in a bounded, smooth domain $\Omega$, then the sequence $(u_i)$ of solutions to the Dirichlet problem for the $p_i(x)$-Laplacian with fixed boundary datum $\varphi$ converges (in a sense to be made precise) to the solution $u_p$ of the Dirichlet problem for the $p(x)$-Laplacian with boundary datum $\varphi$.  A similar result is proved for a decreasing sequence $p_j\searrow p$
	\end{abstract}
	\maketitle
	\newtheorem{question}{Question}
	\newtheorem{property}{Property}
	\newtheorem{theorem}{Theorem}[section]
	\newtheorem{acknowledgement}{Acknowledgement}
	\newtheorem{algorithm}{Algorithm}
	\newtheorem{axiom}{Axiom}[section]
	\newtheorem{case}{Case}
	\newtheorem{claim}{Claim}
	\newtheorem{conclusion}{Conclusion}
	\newtheorem{condition}{Condition}
	\newtheorem{conjecture}{Conjecture}
	\newtheorem{corollary}{Corollary}[section]
	\newtheorem{criterion}{Criterion}
	\newtheorem{definition}{Definition}[section]
	\newtheorem{example}{Example}[section]
	\newtheorem{exercise}{Exercise}
	\newtheorem{lemma}{Lemma}[section]
	\newtheorem{notation}{Notation}
	\newtheorem{problem}{Problem}
	\newtheorem{proposition}{Proposition}[section]
	\newtheorem{remark}{Remark}[section]
	\newtheorem{solution}{Solution}
	\newtheorem{summary}{Summary}

	\thispagestyle{empty}

	\section{Introduction}
	Stability results are of paramount importance in the study of boundary value problems since they address the dependence of the solutions on the data. In the particular case of boundary value problems involving the $p$-Laplacian, there is special interest in the behavior of the solutions $u_p$ with respect to variations of the parameter $p$. For $p$ independent of the space variable, a stability result for the Poisson's problem with vanishing boundary value was studied in \cite{Lindq}. The corresponding generalization to variable $p$ was carried out in \cite{LaMe}. 
	This article is devoted to the discussion of stability for the solutions of the Poisson problem with Dirichlet boundary condition 
	\begin{equation}\label{PP}
		\begin{cases}
			\Delta_{p}(w) = \text{div}\left( |\nabla w|^{p-2} \nabla w \right) = f\in W^{-1,p}(\Omega) & \text{in} \ \Omega, \\
			w|_{\partial \Omega} = \varphi
		\end{cases}
	\end{equation}
	under perturbations of the parameter $p$. In (\ref{PP}), $\Omega$ is a bounded, smooth domain whose boundary and closure are denoted by $\partial\Omega$ and $\overline{\Omega}$, respectively, and $p=p(x)$ is a function on $\Omega$.\\
	The central results presented in this work are theorems \ref{central} and \ref{central1}.
	\begin{theorem}\label{central}
		Let $\Omega\subseteq {\mathbb R}^n$ be a bounded, smooth domain with boundary $\partial \Omega$ and closure $\overline{\Omega}$. Consider a non-decreasing sequence of functions $p_i:\Omega \rightarrow [1,\infty)$, 
		$p_i\in C(\overline{\Omega})$ such that $p_i\rightarrow p$ uniformly in $\Omega$. Fix $\varphi \in W^{1,p}(\Omega)$ and $f\in \left(W_0^{1,(p_1)_-}(\Omega)\right)^{\ast}$. For each $i\in {\mathbb N}$ let $w_i\in W^{1,p_i}(\Omega)$ be the solution to the Poisson's problem (\ref{PP}) corresponding to $p=p_i$. Then
		\begin{itemize}
			\item [$(i)$] There exists $w\in W^{1,p}(\Omega)$ such that the sequence $(w_i)_{i\geq k}$ converges weakly to $w$ in $W^{1,p_k}(\Omega)$, for each $k$,
			\item[$(ii)$] For each $I$,  $\|w_i-w\|_{p_I}\rightarrow 0$ as $i\rightarrow \infty$,\\
			\item[$(iii)$] $\||\nabla (w-w_i)|\|_{p_i}\rightarrow 0$ as $i\rightarrow \infty.$
		\end{itemize}
		Moreover if the limit function $p$ satisfies the log-H\"{o}lder condition (\ref{logholder}), then the limit $w$ belongs to $W^{1,p}(\Omega)$ and it is the unique solution of the problem (\ref{PP}).
	\end{theorem}
	\begin{theorem}\label{central1}
		Let $\Omega$ be as in theorem \ref{central} and let the sequence $(p_i)\in C(\overline{\Omega})$ decrease uniformly to $p$. Let $\varphi \in W^{1,p_1}(\Omega)$ and $f\in \left(W^{1,p}_0(\Omega)\right)^{\ast}$. Denote by $w$ and $w_i$ the unique solution of the problem (\ref{PP}) with exponent $p$ and $p_i$, respectively. Assume $1<p_-$, $p_+<\infty$ and 
		that for some $J\in{\mathbb N}$ it holds that $\int\limits_{\Omega}|\nabla w|^{p_J}dx<\infty.$ Then,
		\begin{equation}\label{strongp}
			\int\limits_{\Omega}|\nabla (w_i-w)|^pdx\rightarrow 0\,\,\text{as}\,\,i\rightarrow \infty.
		\end{equation}
		
	\end{theorem}
	Theorem \ref{central} should be compared to a related result in \cite{MRU}, in which the behavior of the solution $u_p$ of (\ref{PP}) with $f=0$, as $p\rightarrow \infty$ is characterized.
	\section{Notation, terminology and known auxiliary results}
	The main results and definitions in this section have been dealt with in \cite{DHHR,OAP,AOJA,KR,ML}. In the sequel, $\Omega\subset {\mathbb R}^n$ will stand for a bounded, smooth domain with boundary $\partial\Omega$ and $p:\Omega\rightarrow (1,\infty)$ denotes a Borel-measurable function subject to the constraints
	\begin{equation}\label{constraintforp}
		1<p_-=\inf\limits_{\Omega}p(x)\leq \sup\limits_{\Omega}p(x)=p_+<\infty.
	\end{equation}
	For the sake of typographical simplicity, variable exponents will be denoted without specific reference to the spatial variable they depend upon, i.e., $p(x)$ will be written as $p$.\\
	The next definition concerns the variable exponent Lebesgue and Sobolev spaces and their corresponding Luxemburg norm.
	\begin{definition}\cite{DHHR, KR, ML}
		\begin{equation*}
			L^{p}(\Omega)=\left\{f:\int\limits_{\Omega}|\lambda f(x)|^{p}\,dx<\infty,\,\text{for some}\,\,\lambda>0\right\}
		\end{equation*}
		and the Luxemburg norm is defined by
		\begin{equation*}
			\|f\|_{p}=\inf\left\{\lambda>0: \int\limits_{\Omega}\left(|f(x)|/\lambda)\right)^{p}dx\leq 1\right\}.
		\end{equation*}
	\end{definition}
	It is a routine exercise to show that when $p$ is constant on $\Omega$, the above defined spaces coincide with the usual Lebesgue spaces. \\
	If $p\leq q$ are measurable in $\Omega$, the embedding
	$L^q(\Omega)\hookrightarrow L^p(\Omega)$ is continuous, i.e., there exists a positive constant $C(p,q)$ such that
	$\|u\|_p\leq C(p,q)\|u\|_q$ for any $u\in L^q(\Omega)$ \cite{KR}.\\
	The following lemma is elementary and its proof will be omitted \cite{ML}.
	\begin{lemma}\label{normmodular}
		For $p$ satisfying (\ref{constraintforp}) and $u\in L^p(\Omega)$ it holds
		\begin{equation}
			\|u\|_p\leq \max \left\{\left(\int\limits_{\Omega}|u|^pdx\right)^{\frac{1}{p_+}}, \left(\int\limits_{\Omega}|u|^pdx\right)^{\frac{1}{p_-}}\right\}.
		\end{equation}
	\end{lemma}
	\begin{remark}
		{\normalfont It is clear from lemma \ref{normmodular} that 
			\begin{equation}
				\|u\|_p\leq  \left(\int\limits_{\Omega}|u|^pdx\right)^{\frac{1}{p_+}}
			\end{equation}
			if $\int\limits_{\Omega}|u|^pdx\leq 1$ and 
			\begin{equation}
				\|u\|_p\leq  \left(\int\limits_{\Omega}|u|^pdx\right)^{\frac{1}{p_-}}
			\end{equation}
			otherwise.
		}
	\end{remark}
	The following is the variable exponent version of H\"{o}lder inequality and will be used in the sequel. We refer the reader to \cite{KR} for the standard proof.
	\begin{lemma}\label{holder}
		Let $p$ be variable exponents satisfying (\ref{constraintforp}) and $p^{-1}+q^{-1}=1$. Then, for $u\in L^p(\Omega)$ and $v\in L^q(\Omega)$, it holds that
		\begin{equation}
			\int\limits_{\Omega}|uv|dx\leq 2\|u\|_p\|v\|_q.
		\end{equation}
	\end{lemma}
	\begin{definition}\cite{DHHR,KR}
		\begin{equation*}
			W^{1,p}(\Omega)=\left\{f:f\in L^p(\Omega)\,\text{and}\,\, |\nabla f|\in L^p\right\},
		\end{equation*}
		where $|\nabla f|$ stands for the Euclidean norm of $\nabla f$ and the Sobolev norm is defined as
		\begin{equation}\label{sobolevnorm}
			\|f\|_{1,p}=\|f\|_p+\||\nabla f|\|_p.
		\end{equation}
		The Luxemburg norm closure of $C_0^{\infty}(\Omega)$ in $W^{1,p}(\Omega)$ will be denoted as $W^{1,p}_0(\Omega)$.
	\end{definition}
	The following theorem is well known \cite{KR}
	
	\begin{theorem}\label{poincare}
		Let $p\in C(\overline{\Omega})$ satisfy conditions (\ref{constraintforp}).Then the functional
		\begin{equation}
			W^{1,p}_0(\Omega)\ni u\rightarrow \int\limits_{\Omega}|\nabla u|^pdx
		\end{equation}
		is a norm, equivalent to the Sobolev norm (\ref{sobolevnorm}).
	\end{theorem}
	In the sequel, it will be understood that any space $W^{1,q}_0(\Omega)$ is equipped with the norm $u\rightarrow \||\nabla u|\|_q$.\\
	
	\begin{theorem}\label{compactness}
		If $\Omega\subseteq {\mathbb R}^n$ is bounded and $p\in C(\overline{\Omega})$ satisfies the bounds (\ref{constraintforp}), then the embedding
		\begin{equation}
			W^{1,p}_0(\Omega)\hookrightarrow L^p(\Omega)
		\end{equation}
		is compact and the space $W^{1,p}(\Omega)$ is uniformly convex (hence reflexive).
	\end{theorem}
	\begin{proof}
		See \cite{BKM,KR,LaMe}.
	\end{proof}
	The next theorem will play a pivotal role in the proof of theorem \ref{central}.
	\begin{theorem}\label{FZ}
		Assume there exists $M>0$ such that for all $x,y\in \overline{\Omega}$ it holds
		\begin{equation}\label{logholder}
			|p(x)-p(y)|\leq M|\log{|x-y|}|^{-1}.
		\end{equation}
		Then $W^{1,p}_0(\Omega)=W^{1,p}(\Omega)\cap W_0^{1,1}(\Omega).$
	\end{theorem}
	\begin{proof}
		See \cite[Theorem 2.6]{FaZa}, see also \cite{HK}.
	\end{proof}
	
	%{theorem}
	As is customary, the dual space of a Banach space $X$ will be denoted by $X^{\ast}$ and $\langle f,x\rangle$ will stand for the action of $f\in X^{\ast}$ on $x\in X$.\\ It is well known \cite{KR}, that if $\Omega\subset {\mathbb R}^n$ is bounded and $p\leq q$ in $\Omega$ then  $W^{1,q}(\Omega)\hookrightarrow W^{1,p}(\Omega)$  with continuous embedding and thus
	\begin{equation}\label{dualinclusion}
		\left(W^{1,p}(\Omega)\right)^{\ast}\hookrightarrow \left(W^{1,q}(\Omega)\right)^{\ast}
	\end{equation}
	continuously.
	\begin{remark}\label{remarkfunctional}
		{\normalfont Notice that as a consequence of the preceding discussion, if $f\in \left(W_0^{1,p}(\Omega)\right)^{\ast}$ and $(u_j)\subset W_0^{1,q}(\Omega)$ converges weakly to $u\in W_0^{1,q}(\Omega)$, then $\langle f,u_j-u\rangle\rightarrow 0$.  }
	\end{remark}
	The following theorems are particular cases of \cite[Theorem 3.2,Theorem 3.3]{KaMe}, see also \cite{OAP}.
	\begin{theorem}\label{minimizer}
		Let $\Omega\subset {\mathbb R}^n$ be a bounded, smooth domain, $p:\Omega \rightarrow {\mathbb R}$ be a continuous function subject to the constraints (\ref{constraintforp}). Then, for any $\varphi\in W^{1,p}(\Omega)$ and $f\in \left(W^{1,p}_0(\Omega)\right)^{\ast}$ there exists a unique minimizer $u_p\in W^{1,p}_0(\Omega)$ to the Dirichlet integral
		\begin{equation}
			W^{1,p}_0(\Omega)\ni u\rightarrow \int\limits_{\Omega}p^{-1}|\nabla (u-\varphi)|^pdx-\langle f,u\rangle.
		\end{equation}
	\end{theorem}
	Setting $w=\varphi-u$, the following immediate consequence of theorem \ref{minimizer} is obtained.
	\begin{theorem}\label{existencetheorem}
		Let $\Omega\subset {\mathbb R}^n$ be a bounded, smooth domain, $p:\Omega \rightarrow {\mathbb R}$ be a continuous function subject to the constraints (\ref{constraintforp}). Then, for any $f\in \left(W^{1,p}_0(\Omega)\right)^{\ast}$ and any $\varphi\in W^{1,p}(\Omega)$, there exists a unique weak solution $w\in W^{1,p}(\Omega)$ to the Dirichlet problem
		\begin{equation}\label{DP}
			\begin{cases}
				\Delta_{p}w = \text{div}\left( |\nabla w|^{p-2} \nabla w \right) = f & \text{in} \ \Omega, \\
				w|_{\partial \Omega} = \varphi.
			\end{cases}
		\end{equation}
		Specifically, $w$ satisfies the identity
		\begin{equation*}
			-\int\limits_{\Omega}\left|\nabla w\right|^{p-2}\nabla w\nabla h \, dx=\langle f,h\rangle\,\,\,\text{for any}\,\,h \in C^{\infty}_0(\Omega)
		\end{equation*}
		and $w-\varphi \in W^{1,p}_0(\Omega).$
	\end{theorem}
	
	\begin{remark}\label{zero}
		{\normalfont As is apparent from theorems \ref{minimizer} and \ref{existencetheorem}, the null function $0$ is a solution of problem (\ref{DP}) if and only if $\varphi\in W^{1,p}_0(\Omega)$ and $f=0$.}
	\end{remark}
	The following lemma will be used in the sequel.\\
	\begin{lemma} \cite{ELN}\label{epsilonestimate}
		Let $\Omega\subset {\mathbb R}^n$ be a bounded domain and $p$ and $q$ be real-valued-measurable functions on $\Omega$ with $p\leq q\leq p+\epsilon$ in $\Omega$, with $0<\epsilon<1$. Then, for a measurable function $f:\Omega\rightarrow [-\infty,\infty]$ it holds
		\begin{equation}\label{eest}
			\int\limits_{\Omega}|f|^pdx\leq \epsilon|\Omega|+\epsilon^{-\epsilon}\int\limits_{\Omega}|f|^qdx.
		\end{equation}
	\end{lemma}
	\begin{proof}
		\begin{align*}
			\int\limits_{\Omega}|f|^pdx&=\left(\int\limits_{|f|<\epsilon} + \int\limits_{\epsilon<|f|<1}+\int\limits_{|f|>1} \right)|f|^pdx\\ 
			&= I_1+I_2+I_3.
		\end{align*}
		The first integral, $I_1$ is clearly bounded by $\int\limits_{\Omega}\epsilon^pdx\leq \epsilon|\Omega|$. As to $I_2$, observe that
		\begin{align*}
			I_2=\int\limits_{\epsilon<|f|<1}|f|^{p-q}|f|^qdx &\leq \epsilon^{-\epsilon}\int\limits_{\epsilon<|f|<1}|f|^qdx.
		\end{align*}
		In conclusion, 
		\begin{align*}
			\int\limits_{\Omega}|f|^pdx&\leq \epsilon|\Omega|+\epsilon^{-\epsilon}\int\limits_{\epsilon<|f|<1}|f|^qdx+\int\limits_{|f|\geq 1}|f|^qdx\\ &\leq \epsilon|\Omega|+\epsilon^{-\epsilon}\int\limits_{\Omega}|f|^qdx.
		\end{align*}
	\end{proof}
	\begin{corollary}\label{normembedding}
		Under the assumptions of the preceding theorem, for any $f\in L^q(\Omega)$ it holds that
		\begin{equation}\label{normestimate}
			\|g\|_p\leq \left(\epsilon|\Omega|+\epsilon^{-\epsilon}\right)\|g\|_q
		\end{equation}
	\end{corollary}	
	\begin{proof}
		Let $\lambda>0$ be chosen so that $\int\limits_{\Omega}\left|\frac{g}{\lambda}\right|^qdx\leq 1.$	Estimate \ref{eest} holds then for $\frac{g}{\lambda}$ in the place of $f$, yielding
		\begin{equation}
			\left(\epsilon|\Omega|+\epsilon^{-\epsilon}\right)^{-1}	\int\limits_{\Omega}\left|\frac{g}{\lambda}\right|^pdx\leq 1.
		\end{equation}
		Since $\left(\epsilon|\Omega|+\epsilon^{-\epsilon}\right)^{-1}<1$,
		\begin{equation}
			\int\limits_{\Omega}\left|\frac{g}{\left(\epsilon|\Omega|+\epsilon^{-\epsilon}\right)\lambda}\right|^pdx\leq \left(\epsilon|\Omega|+\epsilon^{-\epsilon}\right)^{-1}\int\limits_{\Omega}\left|\frac{g}{\lambda}\right|^pdx\leq 1.
		\end{equation}
		Thus, $\left(\epsilon|\Omega|+\epsilon^{-\epsilon}\right)\lambda\geq\|g\|_p$ and (\ref{normestimate}) follows immediately.
	\end{proof}
	In particular, for $p\leq q\leq p+\epsilon$ and $f\in W^{1,q}(\Omega)$
	\begin{align}\label{sobolevnormestimate}
		\|f\|_{1,p}\leq \left(\epsilon|\Omega|+\epsilon^{-\epsilon}\right)\|f\|_{1,q}
	\end{align}
	and it follows that the norm of the embedding (\ref{dualinclusion}) is at most $\left(\epsilon|\Omega|+\epsilon^{-\epsilon}\right).$
	\begin{corollary}\label{withp}
		Under the assumptions of theorem \ref{epsilonestimate},it holds
		\begin{equation}
			\int\limits_{\Omega}p^{-1}|f|^pdx\leq \epsilon|\Omega|+\epsilon^{-\epsilon}(1+\epsilon)\int\limits_{\Omega}q^{-1}|f|^qdx.
		\end{equation}
	\end{corollary}
	\begin{proof}
		The proof follows from applying estimate (\ref{eest}) to $p^{-\frac{1}{p}}|f|$ and observing that $qp^{-\frac{q}{p}}<(1+\epsilon).$
	\end{proof}
	
	The following elementary, technical lemma will be singled out for use in the sequel.
	\begin{lemma}\label{pwithoutp}
		Let $(p_i)$ $p_i\in C(\overline{\Omega})$, $p_i\rightarrow p$ uniformly in $\Omega$ and $f,f_i \in L^{p_i}(\Omega)$ for each $i\in{\mathbb N}$. Assume that $\int\limits_{\Omega}\frac{|f_i|^{p_i}}{p_i}dx\rightarrow \int\limits_{\Omega}\frac{|f|^{p}}{p}dx$ as $i\rightarrow \infty.$ Then $\int\limits_{\Omega}|f_i|^{p_i}dx\rightarrow \int\limits_{\Omega}|f|^{p}dx$ as $i\rightarrow \infty.$
	\end{lemma}
	
	\begin{proof}
		\begin{align}
			\int\limits_{\Omega}|f_i|^{p_i}dx&=\int\limits_{\Omega}|f_i|^{p_i}\frac{p_i}{p_i}dx=\int\limits_{\Omega}|f_i|^{p_i}\frac{p_i-p}{p_i}dx+\int\limits_{\Omega}|f_i|^{p_i}\frac{p}{p_i}dx \\ \nonumber &=\int\limits_{\Omega}|f_i|^{p_i}\frac{p_i-p}{p_i}dx+\int\limits_{\Omega}\left(\frac{|f_i|^{p_i}}{p_i}-\frac{|f|^{p}}{p}\right)p\,dx+\int\limits_{\Omega}|f|^{p}\,dx.
		\end{align}
	\end{proof}
	The next theorem describes a fundamental geometric property of the functional
	\begin{equation}\label{functional}
		W^{1,p}_0(\Omega) \ni w\rightarrow F(u)=\rho_p(w)=\int\limits_{\Omega}\frac{|\nabla w|^p}{p}dx,
	\end{equation}
	where $|\nabla u|$ stands for the Euclidean norm of $\nabla u$.
	\begin{theorem}\label{uniformconvexity}
		Given a domain $\Omega \subset {\mathbb R}^n$ and a function $p:\Omega\rightarrow {\mathbb R}$ with $p\in C(\overline{\Omega})$ and $p_-=\inf\limits_{\Omega}p(x)>1$, the functional (\ref{functional}) is uniformly convex; this means that for any $\varepsilon: 0<\varepsilon $ there exists $\delta: 0<\delta<1$ such that for any pair $u,v \in W^{1,p}(\Omega)$
		\begin{equation}
			\rho_p\left(\frac{u-v}{2}\right)>\varepsilon\frac{\rho_p(u)+\rho_p(v)}{2}\Rightarrow  \rho_p\left(\frac{u+v}{2}\right)<(1-\delta)\frac{\rho_p(u)+\rho_p(v)}{2}
		\end{equation}
	\end{theorem}
	\begin{proof}
		The proof is rather involved, see \cite{BMB,BKM}. Notice that the boundedness of $p$ is not needed here.
	\end{proof}
	\section{Proof of theorem \ref{central}}
	
	Consider a non-decreasing sequence $(p_i)$ with $p_i\in C(\overline{\Omega})$ and $1< p_i\rightarrow p$ uniformly in $\Omega$. Fix $\varphi \in W^{1,p}(\Omega)$ and $f\in \left(W^{1,(p_1)_-}(\Omega)\right)^{\ast}$; it is assumed that $1<(p_1)_-=\inf\limits_{\Omega}p_1$.\\
	For each $i$, let $u_i$ be the unique minimizer in $W^{1,p_i}_0(\Omega)$ of the functional
	\begin{equation}
		F_i(u)=\int\limits_{\Omega}\frac{|\nabla (u-\varphi)|^{p_i}}{p_i}dx-\langle f,u\rangle,
	\end{equation}
	whose existence is guaranteed by theorem \ref{minimizer}.
	It has been shown in \cite{KaMe} that $w_i=\varphi- u_i$ is the unique solution of the Dirichlet problem (\ref{DP}) with $p=p_i$.\\
	To facilitate the flow of ideas, the proof of theorem \ref{central} will be split into several lemmas.
	\begin{lemma}\label{ujbounded} Let $(u_i)$ be the sequence of minimizers introduced above. Then, for any $J\in{\mathbb N}$, the sequence $(u_i)_{i\geq J}$ is bounded in $W^{1,p_J}_0(\Omega)$.
	\end{lemma}
	\begin{proof}
		Fix $\epsilon:0<\epsilon<e^{-1}$ and a natural number $J$ such that $i\geq J\Rightarrow \|p_i-p\|_{\infty}<\frac{\epsilon}{2}$. Then for $i\geq J$ one has $u_i\in W^{1,p_J}_0(\Omega)$. Assume first that $\int\limits_{\Omega} |\nabla u_i|^{p_J}dx \geq 1$. Then, on account of lemma \ref{normmodular},
		\begin{equation}\label{normJ}
			\||\nabla u_i|\|_{p_J}\leq \left(\int\limits_{\Omega} |\nabla u_i|^{p_J}dx\right)^{\frac{1}{(p_J)_-}}
		\end{equation}
		by virtue of lemma \ref{epsilonestimate}, it follows
		\begin{align}\label{boundforujJ}
			\int\limits_{\Omega} |\nabla u_i|^{p_J}dx &\leq \|p_i-p_J\|_{\infty}|\Omega|+
			\|p_i-p_J\|_{\infty}^{-\|p_i-p_J\|_{\infty}} \int\limits_{\Omega} |\nabla u_i|^{p_i}dx\\ \nonumber & \leq \epsilon|\Omega|+\epsilon^{-\epsilon}\int\limits_{\Omega}
			2^{p_i}\left|\nabla \left(\frac{u_i-\varphi}{2}+\frac{\varphi}{2}\right)\right|^{p_i}dx 
			\\ \nonumber
			&\leq \epsilon|\Omega|+
			\epsilon^{-\epsilon}2^{p_+-1}\left(\int\limits_{\Omega}p_ip_i^{-1}|\nabla (u_i-\varphi)|^{p_i}dx+\int\limits_{\Omega}|\nabla \varphi|^{p_i}dx\right)\\ \nonumber& \leq \epsilon|\Omega|+\epsilon^{-\epsilon}2^{p_+-1}p_+\left(\int\limits_{\Omega}|\nabla( u_i-\varphi)|^{p_i}p_i^{-1}dx-\langle f,u_i\rangle\right)  \\ \nonumber &+ \epsilon^{-\epsilon}2^{p_+-1}\left(p_+\langle f,u_i\rangle+\int\limits_{\Omega}|\nabla \varphi|^{p_i}dx\right).
		\end{align}
		Due to the minimal character of $u_i$, one has
		\begin{align}
			\left(\int\limits_{\Omega}|\nabla( u_i-\varphi)|^{p_i}p_i^{-1}dx-\langle f,u_i\rangle\right)&\leq 
			\int\limits_{\Omega}p_i^{-1}|\nabla \varphi|^{p_i}dx
			\\ \nonumber &\leq \int\limits_{\Omega}|\nabla \varphi|^{p_i}dx
		\end{align}
		and by virtue of lemma \ref{epsilonestimate}, lemma \ref{normmodular} and the choice of $j$, the latter is bounded above by 
		\begin{equation}
			\left(\epsilon|\Omega|+\epsilon^{-\epsilon}\int\limits_{\Omega}|\nabla \varphi|^pdx\right)\leq \left(\epsilon|\Omega|+\epsilon^{-\epsilon}\|\nabla \varphi\|_p^{\alpha}\right),
		\end{equation}
		where $\alpha=p_-$ if $\int\limits_{\Omega}|\nabla \varphi|^pdx\leq 1$ and $\alpha=p_+$ otherwise.\\
		On the other hand, $f\in \left(W^{1,(p_1)_-}_0(\Omega)\right)^{\ast}\subset \left(W^{1,p_J}_0(\Omega)\right)^{\ast}$; also
		\begin{equation}
			\||\nabla u_i|\|_{(p_1)_-}\leq C(J)\||\nabla u_i|\|_{p_J},
		\end{equation}
		for a certain positive constant $C(J)$ independent of $i$. Thus,
		\begin{align}
			\epsilon^{-\epsilon}2^{p_+-1}p_+|\langle f,u_i\rangle|&\leq \epsilon^{-\epsilon}2^{p_+-1}p_+\|f\|_{\left(W^{1,(p_1)_-}_0(\Omega)\right)^{\ast}}\||\nabla u_i|\|_{(p_1)_-}\\ &\leq \epsilon^{-\epsilon}2^{p_+-1}p_+\|f\|_{\left(W^{1,(p_1)_-}_0(\Omega)\right)^{\ast}}C(J)\||\nabla u_i|\|_{p_J}\\ \nonumber &\leq
			\left(\delta^{-1}\epsilon^{-\epsilon}2^{p_+-1}p_+\|f\|_{\left(W^{1,(p_1)_-}_0(\Omega)\right)^{\ast}}\right)^{\frac{(p_J)_-}{(p_J)_--1}}\frac{(p_J)_--1}{(p_J)_-}\\ \nonumber &+
			\frac{1}{(p_J)_-}\left(\delta C(J)\||\nabla u_i|\|_{p_J}\right)^{(p_J)_-}.
		\end{align}
		In particular, for $\delta=\left(\frac{(p_J)_-}{2}\right)^{\frac{1}{(p_J)-}}(C(J))^{-1}$, it follows
		\begin{align}
			\epsilon^{-\epsilon}2^{p_+-1}p_+|\langle f,u_i\rangle|&\leq 
			A+\frac{1}{2}\||\nabla u_i|\|_{p_J}^{(p_J)_-},
		\end{align}
		where 
		\begin{equation*}
			A=\left(\delta^{-1}\epsilon^{-\epsilon}2^{p_+-1}p_+\|f\|_{\left(W^{1,(p_1)_-}_0(\Omega)\right)^{\ast}}\right)^{\frac{(p_J)_-}{(p_J)_--1}}\frac{(p_J)_--1}{(p_J)_-}.
		\end{equation*}
		Thus, (\ref{normJ}), (\ref{boundforujJ}) and the estimates thereafter yield
		\begin{align}\label{final}
			\frac{1}{2}\||\nabla u_i|\|^{(p_J)_-}_{p_J}&\leq \epsilon|\Omega|+\epsilon^{-\epsilon}2^{p_+-1}p_+\left(\epsilon|\Omega|+\epsilon^{-\epsilon}\|\nabla \varphi\|_p^{\alpha}\right)\\
			\nonumber &+A+
			\epsilon^{-\epsilon}2^{p_+-1}\left(\epsilon|\Omega|+\epsilon^{-\epsilon}\|\nabla \varphi\|_p^{\alpha}\right).
		\end{align}
		In conclusion, for any $j\geq J$, 
		\begin{equation}
			\||\nabla u_i|\|_{p_J}\leq \max\left\{1, 2B \right\}^{\frac{1}{(p_J)_-}},
		\end{equation}
		where $B$ is the right hand side of (\ref{final}), which does not depend on $i$.
		
		It follows that the sequence $(u_i)$ is uniformly bounded in $W^{1,p_J}_0(\Omega)$, as claimed.
	\end{proof}
	\begin{lemma}\label{main1}
		In the preceding notation, there is a subsequence of the sequence $(u_i)$ that converges weakly in each $W^{1,p_j}_0(\Omega)$, to a function  $u\in \bigcap\limits_{j=1}^{\infty} W_0^{1,p_j}(\Omega)$.
	\end{lemma}
	\begin{proof}
		Fix $J\in{\mathbb N}$. Theorem \ref{compactness} in conjunction with the theorem of Banach-Alaoglu yields the existence of a function $v_J\in W^{1,p_J}_0(\Omega)$ and a  subsequence $(u_{J+k_j})_{j\geq 1}$ that converges to $v_J$ weakly in $W^{1,p_J}_0(\Omega)$.\\
		Applying the preceding reasoning to the sequence $(u_{J+k_j})_{j\geq m}\subset W^{1,p_{J+k_m}}_0(\Omega)$ for any $m>1$ and denoting its weak limit by $v_{J+m}$ it is immediate that $v_{J+m}=v_J$. Hence, the weak limit does not depend on $J$ and it can be written $v_J=u$. One quickly obtains 
		\begin{equation}\label{uinall}
			u \in \bigcap\limits_{i=1}^{\infty}W_0^{1,p_i}.
		\end{equation}
	\end{proof}
	\begin{remark}\label{subsequence} {\normalfont In the sequel, as is customary, the subsequence $(u_{J+k_j})_{j\geq 1}$ will be denoted by $(u_j)$ and the subsequence $(p_{J+k_j})$ of the sequence $(p_i)$ will be relabeled $(p_j)$.}
	\end{remark}
	The next order of business is to prove that $u\in W^{1,p}(\Omega)$. 
		To that effect, we start with the following assertion.
	\begin{lemma}\label{uinlp}
		The weak limit function $u$ whose existence is proved in lemma \ref{main1} satisfies $\int\limits_{\Omega}|\nabla u|^pdx<\infty$.
	\end{lemma}
	\begin{proof}
		Let $\eta>0$ be arbitrary and fix $v\in W^{1,p}_0(\Omega)$. Choose $0<\epsilon<e^{-1}$ and $0<\delta<1$ such that $\epsilon^{-\epsilon}(1+\epsilon)<(1+\delta)$ and
		\begin{align} \max\left\{3\delta\int\limits_{\Omega}p^{-1}|\nabla (v-\varphi)|^pdx, \epsilon|\Omega|\left(1+\epsilon^{-\epsilon}(1+\epsilon)\right),\delta|\langle f,v\rangle|,\delta|\langle f,u\rangle|\right\} <\frac{\eta}{5};
		\end{align}	
		and $k\in {\mathbb N}$ large enough so that $j\geq k\Rightarrow (1+\delta)|\langle f,u_j-u\rangle|<\frac{\eta}{5}$ and $\|p_j-p||_{\infty}<\frac{\epsilon}{2}$. Because the functional
		\begin{equation}
			W^{1,p_k}_0(\Omega)\ni w\rightarrow \int\limits_{\Omega}p_k^{-1}|\nabla (w-\varphi)|^{p_k}dx
		\end{equation}
		is weakly lower semicontinuous, it follows that 
		\begin{equation}\label{ukintermsofuj}
			\int\limits_{\Omega}p_k^{-1}|\nabla (u-\varphi)|^{p_k}dx\leq \liminf_{j\geq k}\int\limits_{\Omega}p_k^{-1}|\nabla (u_j-\varphi)|^{p_k}dx.
		\end{equation}
		In what follows, let $r(\epsilon)=\epsilon^{-\epsilon}(1+\epsilon)$. Corollary \ref{withp} yields, for $j\geq k$,
		\begin{equation}\label{uk}
			\int\limits_{\Omega}p_k^{-1}|\nabla (u_j-\varphi)|^{p_k}dx\leq 
			\epsilon|\Omega|+r(\epsilon)\int\limits_{\Omega}p_j^{-1}|\nabla (u_j-\varphi)|^{p_j}dx.
		\end{equation} 
		Since $v\in W^{1,p}_0(\Omega)\hookrightarrow W^{1,p_j}_0(\Omega)$,  one has, by virtue of the minimizing property of $u_j$,
		\begin{equation}
			\int\limits_{\Omega}p_j^{-1}|\nabla (u_j-\varphi)|^{p_j}dx-\langle f,u_j\rangle\leq \int\limits_{\Omega}p_j^{-1}|\nabla (v-\varphi)|^{p_j}dx-\langle f,v\rangle
		\end{equation}
		\vspace{-.3in}
		and invoking lemma \ref{epsilonestimate} again one obtains
		\begin{equation}
			\int\limits_{\Omega}p_j^{-1}|\nabla (v-\varphi)|^{p_j}dx\leq \epsilon|\Omega|+r(\epsilon)\int\limits_{\Omega}p^{-1}|\nabla (v-\varphi)|^{p}dx.
		\end{equation} 
		In conclusion, for fixed $k$ as above and $j \geq k$, on account of inequality (\ref{uk}), one has
		\vspace{-.2in}
		\begin{align}
			\int\limits_{\Omega}p_k^{-1}|\nabla (u_j-\varphi)|^{p_k}dx&\leq 
			\epsilon|\Omega|+r(\epsilon)\left(\int\limits_{\Omega}p_j^{-1}|\nabla (u_j-\varphi)|^{p_j}dx-\langle f,u_j\rangle\right)\\ \nonumber &+r(\epsilon)\langle f,u_j\rangle \\ \nonumber &\leq 
			\epsilon|\Omega|+r(\epsilon)\left(\int\limits_{\Omega}p_j^{-1}|\nabla (v-\varphi)|^{p_j}dx-\langle f,v\rangle \right)\\ \nonumber &+r(\epsilon)\langle f,u_j\rangle\\ \nonumber & \leq \epsilon|\Omega|+r(\epsilon)\left(\epsilon |\Omega|+r(\epsilon)\int\limits_{\Omega}p^{-1}|\nabla (v-\varphi)|^pdx-\langle f,v\rangle\right)\\ \nonumber & +r(\epsilon)\langle f,u_j\rangle\\
			\nonumber &\leq
			\epsilon|\Omega|(1+r(\epsilon))+ r(\epsilon)\left(r(\epsilon)\int\limits_{\Omega}p^{-1}|\nabla (v-\varphi)|^pdx-\langle f,v\rangle \right)\\ \nonumber &+r(\epsilon)\langle f,u_j\rangle
			\\ \nonumber &\leq \frac{\eta}{5}+(1+\delta)^2\int\limits_{\Omega}p^{-1}|\nabla (v-\varphi)|^pdx-\langle f,v\rangle-(r(\epsilon)-1)\langle f,v\rangle \\ \nonumber &+r(\epsilon)\langle f,u_j-u\rangle +r(\epsilon)\langle f,u\rangle
			\\ \nonumber & \leq \frac{\eta}{5}+3\delta \int\limits_{\Omega}p^{-1}|\nabla (v-\varphi)|^pdx+\int\limits_{\Omega}p^{-1}|\nabla (v-\varphi)|^pdx-\langle f,v\rangle \\ \nonumber &+\delta |\langle f,v\rangle| +(1+\delta)|\langle f,u_j-u\rangle| +(r(\epsilon)-1)\langle f, u\rangle+\langle f,u\rangle\\
			\nonumber & \leq \int\limits_{\Omega}p^{-1}|\nabla (v-\varphi)|^pdx-\langle f, v\rangle+\langle f, u\rangle + \eta.
		\end{align}
		\vspace{-.2in}
		In all, 
		
		\begin{align}
			\liminf\limits_{j\geq k}\int\limits_{\Omega}p_k^{-1}|\nabla (u_j-\varphi)|^{p_k}dx \leq \int\limits_{\Omega}p^{-1}|\nabla (v-\varphi)|^pdx-\langle f, v\rangle+\langle f, u\rangle + \eta, 
		\end{align}		
		that is, according to (\ref{ukintermsofuj}), for any $\eta>0$, 
		\begin{equation}
			\int\limits_{\Omega}p_k^{-1}|\nabla (u-\varphi)|^{p_k}dx\leq
			\int\limits_{\Omega}p^{-1}|\nabla (v-\varphi)|^pdx-\langle f, v\rangle+\langle f, u\rangle + \eta.
		\end{equation}
		Consequently, 
		\begin{align}
			\int\limits_{\Omega}p^{-1}|\nabla (u-\varphi)|^{p}dx&\leq \liminf_{n\geq k}\int\limits_{\Omega}p_n^{-1}|\nabla (u-\varphi)|^{p_n}dx \\ &\leq 
			\int\limits_{\Omega}p^{-1}|\nabla (v-\varphi)|^{p}dx-\langle f, v\rangle+\langle f, u\rangle + \eta.
		\end{align}
		In particular, the preceding inequality yields
		\begin{equation}
			\int\limits_{\Omega}\frac{|\nabla u|^p}{p}dx\leq 2^{p_+-1}\left(\int\limits_{\Omega}p^{-1}|\nabla (u-\varphi)|^{p}dx+\int\limits_{\Omega}p^{-1}|\nabla \varphi|^{p}dx\right)<\infty,
		\end{equation}
		i.e. $|\nabla u|\in L^{p}(\Omega)$, as claimed.
	\end{proof}
	
	\begin{remark} \label{remark3.2}{\normalfont
			Since $p_j\rightarrow p$ pointwise, Lebesgue's dominated convergence yields
			\begin{equation}\label{intjtoint}
				\int\limits_{\Omega}\frac{|\nabla (u-\varphi)|^{p_j}}{p_j}dx\rightarrow
				\int\limits_{\Omega}\frac{|\nabla (u-\varphi)|^{p}}{p}dx \,\text{as}\,j\rightarrow \infty
			\end{equation}
			and on account of lemma \ref{pwithoutp}, 
			\begin{equation}\label{intjtoint1}
				\int\limits_{\Omega}|\nabla (u-\varphi)|^{p_j}dx\rightarrow
				\int\limits_{\Omega}|\nabla (u-\varphi)|^{p}dx \,\text{as}\,j\rightarrow \infty.
			\end{equation}
		}
	\end{remark}
	The next series of lemmas aims to improve the weak convergence statement $u_j\rightharpoonup u$ as $j\rightarrow \infty.$
	\begin{lemma}\label{UJTOPJ} If $(u_j)$ is the sequence of minimizers defined in lemma \ref{main1}, then
		\begin{equation}\label{ujpjtoup}
			\int\limits_{\Omega}|\nabla (u_j-\varphi)|^{p_j}dx\rightarrow
			\int\limits_{\Omega}|\nabla (u-\varphi)|^{p}dx \,\text{as}\,j\rightarrow \infty.
		\end{equation}		
	\end{lemma}
	\begin{proof}
		Fix $\epsilon, \delta, \eta$ and $k$ as in the proof of Lemma \ref{uinlp}. Observe that due to the minimal character of $u_k$,
		\begin{align}
			\int\limits_{\Omega}\frac{|\nabla (u_k-\varphi)|^{p_k}}{p_k}dx-\langle f, u_k\rangle &\leq \int\limits_{\Omega}\frac{|\nabla (u-\varphi)|^{p_k}}{p_k}dx-\langle f, u\rangle.
		\end{align}
		Taking into consideration the fact that $u_k\rightarrow u$ weakly, it follows from the above that
		\begin{align}\label{suplessthanu}
			\limsup_{k}\int\limits_{\Omega}\frac{\nabla (u_k-\varphi)|^{p_k}}{p_k}dx &\leq \limsup_k \int\limits_{\Omega}\frac{|\nabla (u-\varphi)|^{p_k}}{p_k}dx
		\end{align}
		whereas the fact that $|\nabla u|\in L^p(\Omega)$ coupled with a straightforward application of Lebesgue's dominated convergence theorem, yields
		\begin{align}\label{estu}
			\lim_{k\rightarrow \infty}\int\limits_{\Omega}\frac{|\nabla (u-\varphi)|^{p_k}}{p_k}dx=\int\limits_{\Omega}\frac{|\nabla (u-\varphi)|^{p}}{p}dx.
		\end{align}
		On the other hand, 
		\begin{equation}\label{3.30}
			\int\limits_{\Omega}\frac{|\nabla (u-\varphi)|^{p_k}}{p_k}dx\leq \liminf_{j\geq k}\int\limits_{\Omega}\frac{|\nabla (u_j-\varphi)|^{p_k}}{p_k}dx.
		\end{equation}
		In turn, for $j\geq k$, one has, owing to corollary \ref{withp},
		\begin{align}
			\int\limits_{\Omega}\frac{|\nabla (u_j-\varphi)|^{p_k}}{p_k}dx&\leq 
			\epsilon|\Omega|+\epsilon^{-\epsilon}(1+\epsilon)  \int\limits_{\Omega}\frac{|\nabla (u_j-\varphi)|^{p_j}}{p_j}dx
			\\ \nonumber &\leq \eta+(1+\delta)\int\limits_{\Omega}\frac{|\nabla (u_j-\varphi)|^{p_j}}{p_j}dx.
		\end{align}
		Thus, 
		\begin{align}
			\liminf_{j\geq k}\int\limits_{\Omega}\frac{|\nabla (u_j-\varphi)|^{p_k}}{p_k}dx &\leq \eta+(1+\delta)\liminf_{j\geq k}\int\limits_{\Omega}\frac{|\nabla (u_j-\varphi)|^{p_j}}{p_j}dx
		\end{align}
		and it follows from (\ref{3.30})
		\begin{align}\label{estlimk}
			\lim\limits_{k\rightarrow\infty}\int\limits_{\Omega}\frac{|\nabla (u-\varphi)|^{p_k}}{p_k}dx&=\int\limits_{\Omega}\frac{|\nabla (u-\varphi)|^{p}}{p}dx
			\\ \nonumber &\leq \eta+(1+\delta)\liminf_{j}\int\limits_{\Omega}\frac{|\nabla (u_j-\varphi)|^{p_j}}{p_j}dx.
		\end{align}
		In all, from (\ref{suplessthanu}), (\ref{estu}) and (\ref{estlimk}) one deduces
		\begin{align}\label{3.34}
			\limsup_{k}\int\limits_{\Omega}\frac{\nabla (u_k-\varphi)|^{p_k}}{p_k}dx&\leq \int\limits_{\Omega}\frac{|\nabla (u-\varphi)|^{p}}{p}dx \\ \nonumber& \leq \eta+(1+\delta)\liminf_{j}\int\limits_{\Omega}\frac{|\nabla (u_j-\varphi)|^{p_j}}{p_j}dx, 
		\end{align}
		which in conjunction with lemma \ref{pwithoutp} yields (\ref{ujpjtoup}).
	\end{proof}
	The following, similar result can be proved along the same lines:
	
	\begin{lemma}\label{+ujpjtouP}
		\begin{equation}\label{+ujpjtoup}
			\int\limits_{\Omega}\left|\nabla \left(\frac{u+u_j}{2}-\varphi\right)\right|^{p_j}dx\rightarrow
			\int\limits_{\Omega}|\nabla (u-\varphi)|^{p}dx \,\text{as}\,j\rightarrow \infty.
		\end{equation}
	\end{lemma}
	\begin{proof}
		Since $\frac{u+u_k}{2}\in W^{1,p_k}_0(\Omega)$, the minimality of $u_k$ yields the inequality
		\begin{align} 
			\int\limits_{\Omega}\frac{|\nabla (u_k-\varphi)|^{p_k}}{p_k}dx-\langle f, u_k\rangle &\leq \int\limits_{\Omega}\frac{\left|\nabla \left(\frac{u+u_k}{2}-\varphi\right)\right|^{p_k}}{p_k}dx-\langle f, u\rangle.
		\end{align}
		Thus, owing to (\ref{3.34}), 
		\begin{align} 
			\int\limits_{\Omega}\frac{|\nabla (u-\varphi)|^{p}}{p}dx&=	\limsup_{k}\int\limits_{\Omega}\frac{|\nabla (u_k-\varphi)|^{p_k}}{p_k}dx \leq \limsup_{k}\int\limits_{\Omega}\frac{\left|\nabla \left(\frac{u+u_k}{2}-\varphi\right)\right|^{p_k}}{p_k}dx\\ \nonumber &\leq \limsup_{k}\frac{1}{2}\left( \int\limits_{\Omega}\frac{|\nabla (u-\varphi)|^{p_k}}{p_k}dx+\int\limits_{\Omega}\frac{|\nabla (u_k-\varphi)|^{p_k}}{p_k}dx\right)\\ \nonumber &=\int\limits_{\Omega}\frac{|\nabla (u-\varphi)|^{p}}{p}dx.
		\end{align}
		The latter yields
		\begin{equation}\label{supleqI}
			\limsup_{k}\int\limits_{\Omega}\frac{\left|\nabla \left(\frac{u+u_k}{2}-\varphi\right)\right|^{p_k}}{p_k}dx=\int\limits_{\Omega}\frac{|\nabla (u-\varphi)|^{p}}{p}dx.
		\end{equation}
		It has been shown in the first part of the lemma that it holds
		\begin{align}\label{estu1}
			\int\limits_{\Omega}\frac{|\nabla (u-\varphi)|^{p}}{p}dx &=
			\lim_{k\rightarrow \infty}\int\limits_{\Omega}\frac{|\nabla (u-\varphi)|^{p_k}}{p_k}dx;
		\end{align}
		also
		\begin{align}\label{3.40}
			\int\limits_{\Omega}\frac{|\nabla (u-\varphi)|^{p_k}}{p_k}dx	&\leq \liminf_{j\geq k}\int\limits_{\Omega}\frac{\left|\nabla \left(\frac{u_j+u}{2}-\varphi\right)\right|^{p_k}}{p_k}dx
			\\ \nonumber &\leq \eta+(1+\delta)\liminf_{j}\int\limits_{\Omega}\frac{\left|\nabla \left(\frac{u_j+u}{2}-\varphi\right)\right|^{p_j}}{p_j}dx.
		\end{align}
		From (\ref{estu}) and (\ref{3.40}) one concludes
		\begin{equation}
			\int\limits_{\Omega}\frac{\left|\nabla \left(\frac{u_j+u}{2}-\varphi\right)\right|^{p_j}}{p_j}dx\rightarrow 
			\int\limits_{\Omega}\frac{|\nabla (u-\varphi)|^{p}}{p}dx \,\,\text{as}\,\, k\rightarrow \infty,
		\end{equation}
		which yields lemma \ref{+ujpjtouP} via lemma \ref{pwithoutp}.
	\end{proof}
	The following result yields a stronger convergence result for the sequence of minimizers.
	\begin{theorem}\label{UJPJTOUP}
		Let $(u_j)$ denote the subsequence introduced in remark \ref{subsequence}. Then it holds,
		\begin{equation}\label{ujpjtoUP}
			\lim_{j\rightarrow \infty}\int\limits_{\Omega}\left|\nabla (u-u_j)\right|^{p_j}dx=0.
		\end{equation}
	\end{theorem}
	
	\begin{proof}
		The proof of (\ref{ujpjtoUP}) uses the two preceding statements in conjunction with the uniform convexity property stated in theorem \ref{uniformconvexity}. 
		
		First, observe that by remark \ref{remark3.2} one has
		\begin{equation}\label{up_jtoup}
			\int\limits_{\Omega}p_i^{-1}|\nabla (u-\varphi)|^{p_i}dx\rightarrow 	\int\limits_{\Omega}p^{-1}|\nabla (u-\varphi)|^{p}dx.
		\end{equation}
		On account of (\ref{3.34}), and (\ref{up_jtoup}) it follows that there exists a constant $C(u,\varphi)$, depending only on $u, \varphi$, for which
		\begin{equation}
			\int\limits_{\Omega}p_i^{-1}|\nabla (u-\varphi)|^{p_i}dx+\int\limits_{\Omega}p_i^{-1}|\nabla (u_i-\varphi)|^{p_i}dx\leq C(u,\varphi)=C.
		\end{equation}
		Let $\delta>0$ and assume that for a subsequence $(p_{_{i_k}})$ it holds $\int\limits_{\Omega}\left|\nabla (u_{j_k}-u)\right|^{p_{j_k}}dx > \delta$. Then
		\begin{equation}
			2^{p_+}\int\limits_{\Omega}\left|\frac{\nabla (u_{j_k}-u)}{2}\right|^{p_{j_k}}dx > \delta
		\end{equation}
		and 
		\begin{align}
			\int\limits_{\Omega}\left|\frac{\nabla (u_{j_k}-u)}{2}\right|^{p_{j_k}}dx& > 2^{-p_+}\delta=2^{-p_+}\frac{\delta}{C}C\\ \nonumber &\geq 2^{-p_+}\frac{\delta}{C}\left(\int\limits_{\Omega}p_i^{-1}|\nabla (u-\varphi)|^{p_i}dx+\int\limits_{\Omega}p_i^{-1}|\nabla (u_i-\varphi)|^{p_i}dx\right)\\ \nonumber &
			= 2^{1-p_+}\frac{\delta}{C}\frac{\left(\int\limits_{\Omega}p_i^{-1}|\nabla (u-\varphi)|^{p_i}dx+\int\limits_{\Omega}p_i^{-1}|\nabla (u_i-\varphi)|^{p_i}dx\right)}{2}.
		\end{align}
		On account of theorem \ref{uniformconvexity},it follows that there exists $\eta:0<\eta<1$ for which
		\begin{align}
			\int\limits_{\Omega}p_i^{-1}\left|\nabla \left(\frac{u+u_j}{2}-\varphi\right)\right|^{p_i}dx\leq (1-\eta)\frac{\int\limits_{\Omega}p_i^{-1}|\nabla (u-\varphi)|^{p_i}dx+\int\limits_{\Omega}p_i^{-1}|\nabla (u_i-\varphi)|^{p_i}dx}{2}.
		\end{align}
		Letting $i\rightarrow  \infty$ in the last inequality leads to a contradiction. 
	\end{proof}
	
	The final step is to verify that $u\in W^{1,p}_0(\Omega)$. To this end, it remains to show that $u\in L^p(\Omega)$ and that $u$ can be approximated by $C^{\infty}_0(\Omega)$ functions in the norm $\|\cdot\|_{1,p}$. This is proved in the ensuing proposition:
	
	\begin{proposition}\label{uinw1p0}
		The limit function $u$ obtained in lemma \ref{main1} belongs to $W^{1,p}_0(\Omega).$
	\end{proposition}
	\begin{proof} Observe that since $p_j\nearrow p$, for any $x\in \Omega$ with $p(x)<n$, one has 
		\[ \frac{np_j(x)}{n-p_j(x)}\nearrow \frac{np(x)}{n-p(x)}, \]
		also,
		\[\frac{np(x)}{n-p(x)}-p(x)=\frac{p^2(x)}{n-p(x)}\geq\frac{p_-^2}{n-p_-}>\frac{1}{n-1}. \]
		Consequently, for large enough $j$, 
		\begin{equation}\label{jlargeenough}
			\frac{np_j(x)}{n-p_j(x)}>\frac{np(x)}{n-p(x)}-\frac{1}{n-1}>p(x). 
		\end{equation}
		Let $0<r<\frac{n^2}{p_++n}$. Set $\Omega_1=\{x:1<p(x)<n\}$; $\Omega_2=\{x:n-r<p(x)<p_++r\}$. Then $\Omega_1$ and $\Omega_2$ are open and $\Omega=\Omega_1\cup\Omega_2$; let $\{\chi_1,\chi_2\}$ be a partition of unity subordinated to the cover $\{\Omega_1,\Omega_2\}$. It follows that $u\chi_k\in W^{1,p_j}_0(\Omega_k)$ for $k=1,2$ and all $j\in {\mathbb N}$, whence $u\chi_1\in W^{1,p_j}_0(\Omega_1)\subseteq L^{\frac{np_j}{n-p_j}}(\Omega_1)\subseteq L^p(\Omega)$ for $j$ chosen so that (\ref{jlargeenough}) holds. Similarly, the choice of $r$ yields $n^2-r(p+n)>n^2-r(p_++n)>0$. Thus, $n^2-rn>rp$ and
		\begin{equation}
			u\chi_2\in W^{1,p_j}_0(\Omega_2)\subseteq W^{1,n-r}_0(\Omega_2) \subseteq L^{\frac{n(n-r)}{r}}(\Omega_2)\subseteq L^{p}(\Omega).
		\end{equation}
		In conclusion, $u=u\chi_1+u\chi_2\in L^p(\Omega)$ and it follows that $$u\in W^{1,p}(\Omega)\bigcap\left(\bigcap\limits_{j=1}^{\infty} W^{1,p_j}_0(\Omega)\right).$$ On account of theorem \ref{FZ}, it must hold $u\in W^{1,p}_0(\Omega)$, as claimed.
	\end{proof}
	\begin{corollary}
		The original sequence of minimizers alluded to in lemma \ref{ujbounded} converges weakly to $u\in W^{1,p}_0(\Omega)$ in every $W^{1,p_J}_0(\Omega)$ and theorem \ref{UJPJTOUP} holds for the original sequence $(u_i)$.
	\end{corollary}
	\noindent
	{\bf Proof of theorem \ref{central}}:\\
	The proof of theorem \ref{central} follows by observing that for each $i$, $w_i=\varphi-u_i$ and $w=\varphi -u$.\\
	{\bf Proof of theorem \ref{central1}}:\\
	In this case, the sequence of minimizers $(u_i)$ of the functionals
	\begin{equation}
		W^{1,p_i}_0(\Omega)\ni \rightarrow F_i(v)=\int\limits_{\Omega}\frac{|\nabla (v-\varphi)|^{p_i}}{p_i}dx-\langle f,v\rangle
	\end{equation}
	is in fact bounded in $W^{1,p}_0(\Omega)$. The proof will be sketched, since it follows along the same lines as that of the boundedness of $(u_i)$ in theorem \ref{central}. Due to the assumption on $\varphi$ and on the sequence $(p_i)$ it is immediate from lemma \ref{epsilonestimate} that for $i$ as large as necessary for $\|p-p_i\|_{\infty}<1$, 
	\begin{equation}
		\int\limits_{\Omega}|\nabla \varphi|^{p_i}dx\leq \|p_i-p_1\|_{\infty}|\Omega|+\|p_i-p_1\|_{\infty}^{-\|p_i-p_1\|_{\infty}}\int\limits_{\Omega}|\nabla \varphi|^{p_1}dx.
	\end{equation}
	Also,
	\begin{equation}\label{ujpasujpj}
		\int\limits_{\Omega}|\nabla u_i|^{p}dx\leq \|p_i-p\|_{\infty}|\Omega|+\|p_i-p\|_{\infty}^{-\|p_i-p\|_{\infty}}\int\limits_{\Omega}|\nabla u_i|^{p_i}dx.
	\end{equation}
	Using the same ideas as in lemma \ref{ujbounded}, it can be shown that the sequence $\left(\int\limits_{\Omega}|\nabla u_i|^{p_i}dx\right)$ is bounded. Hence $(u_i)$ is bounded in $W^{1,p}_0(\Omega)$ and therefore that there exists $u\in W^{1,p}_0(\Omega)$ and a subsequence $(u_j)$ of the sequence of minimizers (see remark \ref{subsequence}) such that $u_j\rightharpoonup u$ in $W^{1,p}_0(\Omega)$.\\
	In this setting, one has
	\begin{equation}\label{ujtoup}
		\int\limits_{\Omega}|\nabla (u_j-\varphi)|^pdx\rightarrow \int\limits_{\Omega}|\nabla (u-\varphi)|^pdx\,\,\text{as}\,\,j\rightarrow \infty.
	\end{equation}
	Indeed, owing to the weak lower semicontinuity the minimal character of $u_j$ and the integrability assumption on $u$, one has
	\begin{align}\label{3.56}
		\int\limits_{\Omega}\frac{|\nabla (u-\varphi)|^p}{p}dx &\leq \liminf_{j\geq J}\int\limits_{\Omega}\frac{|\nabla (u_j-\varphi)|^p}{p}dx
	\end{align}
	and for any fixed $0<\epsilon<e^{-1}$ and $J$  so large that $j\geq J$ implies $\|p_j-p\|_{\infty}<\frac{\epsilon}{2}$, one has (corollary \ref{withp})
	\begin{align}\label{3.57}
		\int\limits_{\Omega}\frac{|\nabla (u_j-\varphi)|^p}{p}dx &\leq \epsilon|\Omega|+\epsilon^{-\epsilon}(1+\epsilon)\int\limits_{\Omega}\frac{|\nabla (u_j-\varphi)|^{p_j}}{p_j}dx \\ \nonumber &=\epsilon|\Omega|+\epsilon^{-\epsilon}(1+\epsilon)F_j(u_j)
		+\epsilon^{-\epsilon}(1+\epsilon)\langle f,u_j\rangle 
		\\ \nonumber & \leq \epsilon|\Omega|+\epsilon^{-\epsilon}(1+\epsilon)F_j(u)
		+\epsilon^{-\epsilon}(1+\epsilon)\langle f,u_j\rangle.
	\end{align}
	Taking $\limsup$ in (\ref{3.57}) and observing that since $|\nabla u|\in L^{p_N(\Omega)}$ for some $p_N$
	\begin{equation}
		\limsup_jF_j(u)=\limsup_j\left(\int\limits_{\Omega}\frac{|\nabla (u-\varphi)|^{p_j}}{p_j}dx-\langle f,u\rangle\right)=\int\limits_{\Omega}\frac{|\nabla (u-\varphi)|^{p}}{p}dx-\langle f,u\rangle.
	\end{equation}
	and that $\limsup_j\langle f,u_j\rangle=\langle f,u\rangle$ it follows 
	\begin{equation}
		\limsup_j\int\limits_{\Omega}\frac{|\nabla (u_j-\varphi)|^p}{p}dx\leq \int\limits_{\Omega}\frac{|\nabla (u-\varphi)|^p}{p}dx;
	\end{equation}
	coupled with (\ref{3.56}) and using lemma \ref{pwithoutp}, the latter yiels (\ref{ujtoup}), since $p$ is bounded in $\Omega$. Statement (\ref{strongp}) in theorem \ref{central1} follows from uniform convexity along the same lines as in the proof of theorem \ref{UJPJTOUP}.

\end{document}